\documentclass[a4paper,12pt,ceqn]{article}

\usepackage{t1enc,amssymb,amsbsy,amsthm,amsmath,amsfonts}

\newtheorem*{prop}{Proposition}

\newenvironment{prf}{\textbf{Proof}}{\flushright$\blacksquare$\\}

\newenvironment{sub}{\subsection {}\vspace*{-0.8cm}\hspace*{1.1cm}}{}

\newcommand{\dis}{\displaystyle}

\numberwithin{equation}{section}

\begin{document}
\author{C. Donati-Martin$^{(1)}$, M.Yor$^{(1), (2)}$}

\title{\textbf{Further examples of explicit Krein representations of certain subordinators}}
\date{}
\maketitle

\noindent $^{(1)}$ Laboratoire de Probabilit\'{e}s et Mod\`{e}les
Al\'{e}atoires,\\
\hspace*{0.55cm}Universit\'{e} Paris VI et VII, 4 Place Jussieu - Case 188,\\
\hspace*{0.55cm}F-75252 Paris Cedex 05\\
\\
$^{(2)}$ Institut Universitaire de France \vskip3cm
\noindent\textbf{\underline{Abstract :}} In a previous paper \cite{DY},
we have shown that the gamma subordinators may be represented as
inverse local times of certain diffusions. In the present paper,
we give  such representations for other
subordinators whose Lévy densities are of the form $\dis
\frac{\mathcal{C}}{(\sinh(y))^\gamma}$, $0 < \gamma < 2$, and the more general family obtained from those by  exponential tilting. \newpage
\section{Aim of the paper and summary of \cite{DY}} 
\sub In this paper, we continue the
program started in  \cite{DY}, that is : to represent as many
subordinators $(\mathcal{S}_\ell,\,\ell\geq0)$, i.e : increasing
Lévy processes, started at 0, as possible as inverse local times
$(\tau_\ell,\,\ell\geq0)$ of some particular $\mathbb{R}_+$-valued
diffusion $(X_t)$, such that 0 is regular for itself (relatively
to $X$). More precisely, assume that :
\begin{equation*}
    E[\exp(-\lambda\mathcal{S}_\ell)]=\exp(-\ell \Psi(\lambda))\,,
\end{equation*}
where : $\dis\Psi(\lambda)=\int_0^\infty \nu(dy)(1-e^{-\lambda
y})$, and $\nu(dy)$ - the Lévy measure associated with
$(\mathcal{S}_\ell,\,\ell\geq0)$ - is of the form :
\begin{equation*}
    \nu(dy)=h(y)dy\,,\quad\textrm{with}\quad h(y)=\int_0^\infty
    d\lambda(x)\,e^{-xy}\,,
\end{equation*}
for some positive $\sigma$-finite measure $\lambda(dx)$ on
$\mathbb{R}_+$, then, it is known, as a consequence of Krein's
theory (cf : Knight  \cite{K}, Kotani-Watanabe  \cite{KW}), that there exists a
unique diffusion $(X_t)$ taking values in $\mathbb{R}_+$, such
that its inverse local time at 0, $(\tau_\ell,\,\ell\geq0)$, is
distributed
as $(\mathcal{S}_\ell,\,\ell\geq0)$.\\
Finding $X$ when $\nu$ is given is called (here) \emph{Krein
representation problem}.\\
In our paper  \cite{DY}, we could fill in the following
\begin{center}\underline{\textbf{Table 1}}\\~\\
$\begin{array}{c|c|c}
 && \\
  h(y) & \textrm{Generator of }(X_t) &  \textrm{Distribution }\\
  &&\\
  \hline
  &&\\
  \frac{\mathcal{C}}{y^{\alpha+1}} & \mathcal{L}_{-\alpha}=\frac{1}{2}\frac{d^2}{dx^2}+\frac{\delta-1}{2x}\frac{d}{dx}\,;\quad\delta=2(1-\alpha)& P^\delta
  \\(0<\alpha<1) && \\
  \hline
  && \\
  \frac{\mathcal{C}}{y^{\alpha+1}}\,e^{-\mu y} & \mathcal{L}_{-\alpha}^{\mu\downarrow}=\mathcal{L}_{-\alpha}+\sqrt{2\mu}\,\frac{\widehat{K}'_\alpha(\sqrt{2\mu}x)}{\widehat{K}_\alpha(\sqrt{2\mu}x)}\frac{d}{dx}& P^{\delta; \mu\downarrow} \\
  (0<\alpha<1;\mu>0)&\textrm{where }\widehat{K}_\alpha(y)=y^\alpha K_\alpha(y),\,y>0 & \\ 
  &&\\
  \hline
  && \\
  \frac{\mathcal{C}}{y}\,e^{-\mu y} & \mathcal{L}_{0}^{\mu\downarrow}=\frac{1}{2}\frac{d^2}{dx^2}+\left(\frac{1}{2x}+\sqrt{2\mu}\frac{K'_0}{K_0}\,(\sqrt{2\mu}x)\right)\frac{d}{dx} & P^{2; \mu \downarrow}
  \\(\mu>0)&& \\
\end{array}$
\end{center}
\subsection{Details of Table 1} In fact, the result for the first
row goes back at least to Molchanov-Ostrovski  \cite{MO}, the result for
the second row is deduced from that in the first row with the help
of the following discussion, which relates Esscher transforms (of
subordinators) to Girsanov transforms (of diffusions).\\
If $(X_t)_{t\geq0}$ is a $\mathbb{R}_+$-valued diffusion, whose
inverse local time at 0 :
\begin{equation*}
\tau_\ell =\inf\{t:L_t>\ell\},\ell\geq0\;,
\end{equation*}
admits Lévy measure $\nu(dy)$, and Lévy exponent
$(\Psi(\theta),\theta\geq0)$, and if one defines :
\begin{equation*}
    \varphi_{\theta\downarrow}(x)=E_x[\exp(-\theta T_0(X))]\,,
\end{equation*}
then, there is another diffusion, which we shall denote by
$(X_t^{\theta\downarrow},t\geq0)$, with laws
$(P_x^{\theta\downarrow},x\geq0)$, such that :
\begin{equation}\label{eq1.1}
    P_{x|_{\mathcal{F}_t}}^{\theta\downarrow}=\frac{\varphi_{\theta\downarrow}(X_t)}{\varphi_{\theta\downarrow}(x)}\exp(\Psi(\theta)L_t-\theta
    t)\textrm{\huge.}\,P_{x|_{\mathcal{F}_t}}
\end{equation}
whose inverse local time $(\tau_\ell,\ell\geq0)$ under
$P_0^{\theta\downarrow}$ satisfies :
\begin{equation*}
    E_0^{\theta\downarrow}(\exp-\lambda\tau_\ell)=\exp(-\ell(\Psi(\lambda+\theta)-\Psi(\theta)))
\end{equation*}
i.e : this inverse local time is the $\theta$-Esscher transform of
$(\tau_\ell,\ell\geq0)$ under \mbox{$P_0$ :} its Lévy measure
(under $P_0^{\theta\downarrow}$) is :
$e^{-\theta y}\nu(dy)$.\\
It is also noteworthy that, under some adequate restriction of their domains, the
infinitesimal generators $\mathcal{L}^{\theta\downarrow}$ and
$\mathcal{L}$ are related by :
\begin{equation*}
    \mathcal{L}^{\theta\downarrow}=\mathcal{L}+\frac{d}{dx}\left(\log(\varphi_{\theta\downarrow}(x))\right)\textrm{\huge.}\,\frac{d}{dx}
\end{equation*}
Finally, the result for the third row was deduced by letting
$\alpha\to0$ in the second row, while taking care of the choice\footnote{As is well-known, the local time in a standard Markovian set up, at a given level, is unique up to a multiplicative constant, which for our studies, needs to be chosen carefully.} of
the local times made for
$\mathcal{L}_{-\alpha\,,\,\mu\downarrow}$. (A compendium of choices of local times for Bessel-like diffusions is made in \cite{DRVY}).
\sub In the present
paper, we wish to complete the preceding Table 1, by considering
the 3 parameter family of Lévy measures on $\mathbb{R}_+$ :
\begin{equation}\label{eq1.2}
    \nu_{\mu,\,\alpha,\,k}(dy)=\mathcal{C} \left(\frac{\mu}{\sinh(\mu
    y)}\right)^{\alpha+1}\exp(\mu ky) \,dy
\end{equation}
(The "true" parameters are : $\mu>0$, $k$, and $\alpha$; as
before, $\mathcal{C}$ is simply there to ensure an additional
degree of freedom, if necessary).\\
In order that $\nu_{\mu,\,\alpha,\,k}(dy)$ be a Lévy measure, i.e
: it must satisfy
\begin{equation*}
\int_0^\infty
(x\wedge1)\,\nu_{\mu,\,\alpha,\,k}\,(dx)<\infty\,,\quad\textrm{we
need :}\quad0\leq\alpha<1\,;\;k<1+\alpha\,.
\end{equation*}
We now recall that, from Pitman-Yor  \cite{PY3} formulae (16), p. 276),
if $Q_z^{\delta,\,\mu}$, for $0<\delta\equiv2(1-\alpha)<2$, and
$\mu>0$, denotes the distribution of the squared radial
Ornstein-Uhlenbeck process, with "dimension" $\delta$, and
parameter $\mu$, started from $z$, i.e : the solution of :
\begin{equation*}
    dZ_t=2\sqrt{Z_t}dB_t+(\delta-2\mu
    Z_t)dt\,;\qquad Z_t\geq0\,,\;Z_0=z\,,
\end{equation*}
then, under $Q_0^{\delta,\mu}$, the inverse local time
$(\tau_\ell,\ell\geq0)$ admits as its Lévy \mbox{measure :}
\begin{equation}\label{eq1.3}
    \mathcal{C} \left(\frac{\mu}{\sinh(\mu y)}\right)^{\alpha+1}\exp(\mu\frac{\delta}{2}y) \,dy
\end{equation}
which is a particular case of \eqref{eq1.2}, with $\dis k=\frac{\delta}{2}=(1-\alpha)$.\\
In the next section, we shall show, essentially with the help of
the recipe \eqref{eq1.1}, how to construct a diffusion, indexed by the 3
parameters $(\alpha,\mu,k)$, which solves \underline{the
Krein representation problem for} $\nu_{\mu,\,\alpha,\,k}$.
\sub Some
among the new diffusions we are finding as solutions of Krein's
problem are related to the diffusions we found in  \cite{DY} by time
changing. We first discovered this relationship by applying the
analytical identity :
\begin{equation}\label{eq1.4}
    W_{0,\,\beta}(z)=\sqrt{\frac{z}{\pi}}\;K_\beta\;\left(\frac{z}{2}\right)
\end{equation}
between $W_{0,\,\beta}\;$, a Whittaker function with parameters
$(0,\,\beta)$, and $K_\beta$ (see Appendix). Thus, a part of our
present discussion may be considered as giving a probabilistic
interpretation to \eqref{eq1.4}.\\
We also develop a similar discussion for the analytical  identity
$$  M_{0,\,\beta}(z)= 4^{\beta} \Gamma(\beta +1) \sqrt{z}\;I_{\beta}\left(\frac{z}{2}\right).$$
\section{Solving Krein's problem for $\nu_{\mu,\,\alpha,\,k}$}
\sub We take up the notation in \eqref{eq1.3}; in fact, it is more
convenient to consider the family of radial Ornstein-Uhlenbeck
processes (and not their squares), which we shall denote as
$(R^{\delta,\,\mu}(t),\,t\geq0)$ and their laws
$(P^{\delta,\,\mu}_r,$\\$r\geq0)$. It will be helpful, for the
sequel, to have the following formula at hand, for the
infinitesimal generator $\mathcal{L}_{-\alpha,\,\mu}$ of
$R^{\delta,\,\mu}$ :
\begin{equation}\label{eq2.1}
    \mathcal{L}_{-\alpha,\,\mu}=\frac{1}{2}\;\frac{d^2}{dx^2}\;+\left(\frac{\delta-1}{2x}-\mu
    x\right)\;\frac{d}{dx}
\end{equation}
It is well-known (see, e.g., Pitman-Yor  \cite{PY1}, p. 454, formula
(6.b)) that there is the relationship :
\begin{equation}\label{eq2.2}
    R^{\delta,\,\mu}(t)=e^{-\mu t}R^\delta \left(\frac{e^{2\mu
    t}-1}{2\mu}\right)\,,\qquad t\geq0\,,
\end{equation}
where, on the RHS, $(R^\delta(u),\,u\geq0)$ denotes a
$\delta$-dimensional Bessel process. Thus, we obtain :
\begin{equation}\label{eq2.3}
    \left(\frac{e^{2\mu
    T_0}-1}{2\mu}\,;\;P^{\delta,\,\mu}_x\right)\underset{(a)}{\overset{(law)}{=}}(T_0\,;\;P_x^\delta)\underset{(b)}{\overset{(law)}{=}}\frac{x^2}{2\gamma_\alpha}
\end{equation}
where, on the RHS, $\gamma_\alpha$ denotes a gamma variable with
parameter $\alpha$. [(a) follows from \eqref{eq2.2}, while (b) is
well-known, and goes back to Getoor  \cite{G}, see, e.g., Yor  \cite{Y}, for
some variants...].\\
We thus deduce the following formula from \eqref{eq2.3}, with the help of
elementary computations :
\begin{eqnarray}
 \lefteqn{ E^{\delta,\,\mu}_x(\exp(-\theta T_0))=\frac{1}{\Gamma(\alpha)(\mu x^2)^{\frac{\theta}{2\mu}}}\int_0^\infty \frac{t^{\alpha-1+\frac{\theta}{2\mu}}e^{-t}}{(1+\frac{t}{\mu x^2})^{\frac{\theta}{2\mu}}} \ dt} \nonumber\\
  & = \frac{\Gamma(\alpha+\frac{\theta}{2\mu})}{\Gamma(\alpha)}(\mu
  x^2)^{\frac{\alpha-1}{2}}e^{\mu\frac{x^2}{2}}W_{\frac{(1-\alpha)-\frac{\theta}{\mu}}{2}\,,\;\frac{\alpha}{2}}(\mu x^2) \quad  \mbox{ for } -2\alpha \mu < \theta & \label{eqlaplace}
\end{eqnarray}
where $W_{a,\,b}$ denotes the Whittaker function, with parameters
$(a,b)$.
\sub We now write :
\begin{eqnarray*}
     \nu_{\mu,\,\alpha,\,k}(dy)&=&\mathcal{C}\,\left(\frac{\mu}{\sinh(\mu y)}\right)^{\alpha+1}\exp(\mu ky)dy\\
     &\equiv&\mathcal{C}\,\left(\frac{\mu}{\sinh(\mu y)}\right)^{\alpha+1}\exp(\mu\frac{\delta}{2}y)\exp(-\theta
     y)dy\,,
\end{eqnarray*}
where : $\dis \theta=\mu\left(\frac{\delta}{2}-k\right), \ k < 1+ \alpha$.\\
According to the preceding computation, we now find that the
diffusion with infinitesimal generator :
\begin{equation}\label{eq2.4}
    \mathcal{L}_{-\alpha,\,\mu}^{\theta\downarrow}\equiv\mathcal{L}_{-\alpha,\,\mu}+\frac{d}{dx}\log\left\{(\mu x^2)^{\frac{\alpha-1}{2}}e^{\mu\frac{x^2}{2}}W_{\frac{(1-\alpha)-\frac{\theta}{\mu}}{2}\,,\;\frac{\alpha}{2}}(\mu x^2)\right\}\textrm{\huge.}\frac{d}{dx}
\end{equation}
solves Krein's representation problem for
$\nu_{\mu,\,\alpha,\,k}$. (We note in fact that : $\dis
\frac{(1-\alpha)-\frac{\theta}{\mu}}{2}=\frac{k}{2}$, so that :
$\dis
W_{\frac{(1-\alpha)-\frac{\theta}{\mu}}{2},\,\frac{\alpha}{2}}(\xi)\equiv
W_{\frac{k}{2},\,\frac{\alpha}{2}}(\xi)$).\\
The case where $k=0$ is particularly interesting, since, on one
hand :
\begin{equation*}
     \nu_{\mu,\,\alpha,\,0}(dy)=\mathcal{C}\,\left(\frac{\mu}{\sinh(\mu y)}\right)^{\alpha+1}dy\,,
\end{equation*}
and on the other hand (see Appendix) :
\begin{equation}\label{eq2.5}
   W_{0,\,\frac{\alpha}{2}}(\xi)=\sqrt{\frac{\xi}{\pi}}\;K_{\frac{\alpha}{2}}\left(\frac{\xi}{2}\right)
\end{equation}
so that the diffusion which solves Krein's representation problem
for $\nu_{\mu,\,\alpha,\,0}(dy)$ is the solution to :
\begin{equation}\label{eq2.6}
    dX_t=dB_t+\left[ \frac{\delta-1}{2X_t}+\mu
    X_t\left(\frac{\widehat{K}'_{\frac{\alpha}{2}}}{\widehat{K}_{\frac{\alpha}{2}}}\right)\left(\mu\frac{X_t^2}{2}\right) \right] dt
\end{equation}
Here, we need to give some details about this computation :
\begin{enumerate}
    \item [a)] We deduce from formula \eqref{eq2.4}, in the particular
    case $k=0$, i.e : $\dis \theta=\frac{\delta\mu}{2}$ with the
    help of formula \eqref{eq2.5}, that :
    \begin{equation*}
        \mathcal{L}_{-\alpha,\,\mu}^{\theta\downarrow}=\mathcal{L}_{-\alpha,\,\mu}+\left(\frac{\alpha}{x}+\mu x+\mu x\frac{K'_{\frac{\alpha}{2}}\left(\mu\frac{x^2}{2}\right)}{K_{\frac{\alpha}{2}}\left(\mu\frac{x^2}{2}\right)}\right)\textrm{\huge.}\,\frac{d}{dx}
    \end{equation*}
    \item [b)] Now, trivially :
    \begin{equation*}
        \mathcal{L}_{-\alpha,\,\mu}+\mu x\,\frac{d}{dx}=\mathcal{L}_{-\alpha}\;,
    \end{equation*}
    and, equally simply :
    \begin{equation*}
        \frac{\alpha}{x}+\mu x\frac{K'_{\frac{\alpha}{2}}\left(\mu\frac{x^2}{2}\right)}{K_{\frac{\alpha}{2}}\left(\mu\frac{x^2}{2}\right)}=\mu x\frac{\widehat{K}'_{\frac{\alpha}{2}}}{\widehat{K}_{\frac{\alpha}{2}}}\left(\mu\frac{x^2}{2}\right)
    \end{equation*}
    which translates into the stochastic differential equation
    form of formula \eqref{eq2.6}.
\end{enumerate}
\newpage
As an introduction to the next discussion, we write down
\begin{center}\underline{\textbf{Table 2}}\\~\\
$\begin{array}{c|c|c}
&&\\
   h(y) & \textrm{Generator of }(X_t) & \textrm{Distribution }\\
   &&\\
  \hline
  &&\\
  \frac{\mathcal{C}}{y^{\alpha+1}}\,e^{-\mu y} & \mathcal{L}_{-\alpha}^{\mu\downarrow}=\mathcal{L}_{-\alpha}+\sqrt{2\mu}\,\frac{\widehat{K}'_\alpha}{\widehat{K}_\alpha}(\sqrt{2\mu}x)\frac{d}{dx} & P^{\delta; \mu \downarrow} \\
  (0<\alpha<1;\mu>0)& & (\delta = 2(1-\alpha))\\
  &&\\
  \hline
  &&\\
   \mathcal{C}\left(\frac{\mu}{\sinh(\mu y)}\right)^{\alpha+1} e^{\mu\frac{\delta}{2}y} & \mathcal{L}_{-\alpha,\,\mu}\equiv\mathcal{L}_{-\alpha}- \mu x \frac{d}{dx} & P^{\delta, \mu} \\
 0<  \delta = 2(1-\alpha)<2; \mu >0&&\\
    &&\\
   \hline
   &&\\
  \mathcal{C}\left(\frac{\mu}{\sinh(\mu y)}\right)^{\alpha+1} & \mathcal{L}_{-\alpha,\,\mu}^{\theta\downarrow}\equiv\mathcal{L}_{-\alpha}+\mu
x\,\frac{\widehat{K}'_{\frac{\alpha}{2}}}{\widehat{K}_{\frac{\alpha}{2}}}\left(\mu\frac{x^2}{2}\right)\frac{d}{dx} & P^{\delta, \mu; \frac{\delta \mu}{2} \downarrow}  \\
  &\left(\theta=\frac{\delta\mu}{2}\right)&\\
  &&\\
\end{array}$
\end{center}
The first row is simply taken from Table 1 (second row there).\\
As said above, the second row follows from Pitman-Yor \cite{PY3}.
In the third row, we have written
$\mathcal{L}_{-\alpha,\,\mu}^{\theta\downarrow}$ for the
infinitesimal generator of the process which is defined as : the
radial Ornstein-Uhlenbeck process, with dimension
$\delta=2(1-\alpha)$, and drift parameter $(-\mu)$, pushed
downwards with parameter $\dis\theta=\frac{\delta\mu}{2}$. That
this infinitesimal generator may be expressed in terms of $\dis
\widehat{K}_{\frac{\alpha}{2}}$ will be discussed after (2.7).
\sub We shall now prove a remarkable relationship between the two
families of diffusions whose infinitesimal generators are found on
the RHS of Table 2. This relationship explains precisely why (Row
1) may be deduced from (Row 2), and vice-versa.
\begin{prop}~\\
The following relationship holds with : $\dis
\theta=\frac{\delta\mu}{2}$ :
\begin{equation}\label{eq2.7}
    {X^2}_{-\alpha,\,\mu;\,\theta\downarrow}(t)=X_{-\frac{\alpha}{2};\,\left(\frac{\mu^2}{8}\right)\downarrow}\left(4\int_0^t X^2_{-\alpha,\,\mu;\,\theta\downarrow}(u) \ du\right)
\end{equation}
meaning that : starting from $X\equiv
X_{-\alpha,\,\mu;\,\theta\downarrow}$ on the LHS, there exists\\
$\dis
\left(X_{-\frac{\alpha}{2};\,\frac{\mu^2}{8}}(u),\,u\geq0\right)$
such that the relationship  \eqref{eq2.7} holds.
\end{prop}
\noindent \textbf{Comment about our notation}\\
In formula \eqref{eq2.7}, and possibly several times below, we have
written $X_{i;\,\theta\downarrow}$, etc... instead of
$X_i^{\theta\downarrow}$, for some index $i$. It seemed more
appropriate here, because of the power 2 on the left-side of
(2.7).\\
There should be no confusion between the different diffusions
$X_{i,\,\theta}$ and
$X_{i;\,\theta\downarrow}$.\\\\
\begin{prf}~\\
We start from the stochastic differential equation satisfied by $\dis
(X_{-\alpha,\,\mu;\,\theta\downarrow}(t),$\\$t\geq0)$ as described
(implicitly) in Row 2 of Table 2.\\
Then, taking squares, we obtain :
\begin{equation*}
    X_t^2=x^2+2\int_0^tX_sdB_s+\delta t+2\int_0^t(\mu X_s^2)\frac{\widehat{K}'_{\frac{\alpha}{2}}}{\widehat{K}_{\frac{\alpha}{2}}}\left(\mu\frac{X_s^2}{2}\right)\ ds
\end{equation*}
We now define $(Y_u \equiv Y(u),\,u\geq0)$ via :
\begin{equation*}
    X_t^2=Y\left(4\int_0^t X_s^2 ds\right)\;,\quad(t\geq0)
\end{equation*}
and find that $Y$ satisfies :
\begin{equation*}
    Y_u=x^2+\beta_u+\left(\frac{\delta}{4}\right)\int_0^u\frac{ds}{Y_s}+\frac{1}{2}\int_0^u \mu\;\frac{\widehat{K}'_{\frac{\alpha}{2}}}{\widehat{K}_{\frac{\alpha}{2}}}\left(\frac{\mu}{2}Y_s\right)\ ds
\end{equation*}
since $\dis \frac{\delta}{4}=\frac{\widehat{\delta}-1}{2}$, with
$\dis \widehat{\delta}=2-\alpha=2(1-\frac{\alpha}{2})$, we find
that $(Y_u,\,u\geq0)$ is precisely the diffusion with
infinitesimal generator $\mathcal{L}_{-\frac{\alpha}{2}}^{\nu\downarrow}$,
with $\dis \sqrt{2\nu}=\frac{\mu}{2}$, i.e
$\dis\nu=\frac{\mu^2}{8}$.
\end{prf}
\noindent We now remark that the proof we have just given for the
Proposition relies upon the identification of the
infinitesimal generator of the diffusion $\dis
X_{-\alpha,\,\mu;\,\left(\frac{\delta\mu}{2}\right)\downarrow}$
as given in Table 2; this identification was obtained
from an analytical identity between $W_{0,\textrm{\huge.}}$
and $K_{\textrm{\huge.}}$. (see formula (2.5)).\\
We now explain and prove the Proposition without relying on such
identities, but rather on absolute continuity relationships
between the different laws involved.\\
We now find it a little more convenient to refer to the laws
$\{Q_z^{\delta,\,\mu}\}$ and the main absolute continuity result
we need is :
\begin{equation}\label{eq2.8}
    Q_{z|_{\mathcal{F}_t}}^{\delta,\,\mu}=\exp\left(-\frac{\mu}{2}(Z_t-\delta
    t-z)-\frac{\mu^2}{2}\int_0^tZ_sds\right)\textrm{\huge.}\,Q_{z|_{\mathcal{F}_t}}^\delta\;.
\end{equation}
Here, $(Z_t, t \geq 0)$ denotes the coordinate process on the canonical space $C(\mathbb{R}_+,\mathbb{R}_+)$. \\
We now combine this relation \eqref{eq2.8} with that of the "push downwards"
with parameter $\theta$, so that, with notations which we shall
explain after writing the formula :
\begin{eqnarray}\label{eq2.9}
    &&Q_{z|_{\mathcal{F}_t}}^{\delta,\,\mu;\,\theta\downarrow}\\
    &=&\frac{\varphi_{\theta\downarrow}(Z_t)}{\varphi_{\theta\downarrow}(z)}\exp(\Psi(\theta)L_t-\theta t)\exp\left(-\frac{\mu}{2}(Z_t-\delta t-z)-\frac{\mu^2}{2}\int_0^tdsZ_s\right)\textrm{\huge.}\,Q_{z|_{\mathcal{F}_t}}^\delta\nonumber
\end{eqnarray}
and we note that for precisely : $\dis\theta=\frac{\delta\mu}{2}$,
this relation simplifies as :
\begin{equation}\label{eq2.10}
    Q_{z|_{\mathcal{F}_t}}^{\delta,\,\mu;\,\left(\frac{\delta\mu}{2}\right)\downarrow}=\frac{\exp\left(-\frac{\mu}{2}Z_t\right)\varphi_{\theta\downarrow}(Z_t)}{\exp\left(-\frac{\mu}{2}z\right)\varphi_{\theta\downarrow}(z)}\exp\left(\Psi(\theta)L_t-\frac{\mu^2}{2}\int_0^tdsZ_s\right)\textrm{\huge.}\,Q_{z|_{\mathcal{F}_t}}^\delta
\end{equation}
(The due explanation of the formula \eqref{eq2.9} is that we have combined
the "push-downwards" formula \eqref{eq1.1}, relative to
$\{Q_z^{\delta,\,\mu}\}$, - i.e. the function
$\varphi_{\theta\downarrow}$ and $\Psi(\theta)$ are relative to
that diffusion - with the preceding formula \eqref{eq2.8}). From now on, we
keep : $\dis\theta=\frac{\delta\mu}{2}$. 
\sub We now consider what
becomes of formula \eqref{eq2.10}, once we time change both sides with the
inverse of $\dis\left(4\int_0^tZ_udu,\,t\geq0\right)$, so that, by a slight abuse of notation, the
process of reference is now $(\widehat{Z}(h),\,h\geq0)$, with
$\widehat{Z}$ defined by~:
\begin{equation}\label{eq2.11}
    \framebox{$\dis Z_t=\widehat{Z}\left(4\int_0^tZ_udu\right)$}
\end{equation}
Thus, we obtain :
\begin{equation}\label{timechange}
\widehat{Q}_{z|_{\mathcal{\widehat{F}}_u}}^{\delta,\,\mu;\,\left(\frac{\delta\mu}{2}\right)\downarrow}=
\frac{\exp\left(-\frac{\mu}{2}\widehat{Z}_u\right)\varphi_{\theta\downarrow}(\widehat{Z}_u)}{\left(e^{-\frac{\mu}{2}z}\varphi_{\theta\downarrow}(z)\right)}\exp\left(\Psi(\theta)\widehat{L}_u-\frac{\mu^2}{8}u\right)\textrm{\huge.}\,\widehat{Q}_{z|_{\mathcal{\widehat{F}}_u}}^\delta
\end{equation}
From the well-known property of time change for Bessel processes (see \cite{RY}, Chap.XI, Prop.1.11), $\widehat{Q}_{z}^\delta$ is the distribution of a Bessel process of index $\alpha/2$, i.e. of dimension $\hat{\delta} = 2-\alpha$, that is 
$\widehat{Q}_{z}^\delta = P_z^{\hat{\delta}}$.
Again, with obvious notation, the  right-hand side of \eqref{timechange} may be written
:
\begin{equation*}
    \frac{\widehat{\varphi}_{\frac{\mu^2}{8}\downarrow}(\widehat{Z}_u)}{\widehat{\varphi}_{\frac{\mu^2}{8}\downarrow}(z)}\exp\left(\widehat{\Psi}\left(\frac{\mu^2}{8}\right)\widehat{L}_u-\frac{\mu^2u}{8}\right)\textrm{\huge.}\,\widehat{Q}_{z|_{\mathcal{\widehat{F}}_u}}^\delta\;,
\end{equation*}
and we discover that :
\begin{equation}\label{eq2.12}
\left\{%
\begin{array}{l}
    \widehat{\varphi}_{\frac{\mu^2}{8}\downarrow}(z)=e^{-\frac{\mu}{2}z}\varphi_{\theta\downarrow}(z) \\
    \\
    \widehat{\Psi}\left(\frac{\mu^2}{8}\right)=\Psi(\theta) \\
\end{array}%
\right.
\end{equation}
and $ \widehat{Q}_{z}^{\delta,\,\mu;\,\left(\frac{\delta\mu}{2}\right)\downarrow} = P_z^{\hat{\delta}; \frac{\mu^2}{8} \downarrow}$. 

\noindent
Again, let us explain, very much in the same spirit, e.g : the
first relation :
$\dis\widehat{\varphi}_{\frac{\mu^2}{8}\downarrow}(z)=e^{-\frac{\mu}{2}z}\varphi_{\theta\downarrow}(z)$
in \eqref{eq2.12}.\\
This translates  as :
\begin{equation}\label{eq2.13}
    \widehat{E}_z\left(\exp-\frac{\mu^2}{8}T_0(\widehat{Z})\right)=e^{-\frac{\mu}{2}z}E_z\left(e^{-\frac{\delta\mu}{2}T_0(Z)}\right)
\end{equation}
where $\widehat{Z}$ simply denotes a BES process with dimension
$\widehat{\delta}$, and $Z$ a process with law
$Q_z^{\delta,\,\mu}$. This may be well understood by considering
the absolute continuity relationship \eqref{eq2.8}, when we replace $t$ by
$T_0(Z)$. Then, it follows from that relationship that :
\begin{eqnarray*}
    Q_z^{\delta,\,-\mu}\left(\exp\left(-\frac{\delta\mu}{2}T_0(Z)\right)\right)&=&e^{\frac{\mu z}{2}}Q_z^\delta\left(e^{-\frac{\mu^2}{2}\int_0^{T_0}ds\,Z_s}\right)\\
    &=&e^{\frac{\mu z}{2}}Q_z^{\widehat{\delta}}\left(e^{-\frac{\mu^2}{8}T_0(\widehat{Z})}\right)\;,
\end{eqnarray*}
which is precisely \eqref{eq2.13}. \\
Now, it is well known (see \cite{G}, \cite{Ke}, \cite{PY2}) that the Laplace transform of $T_0$, under the distribution $P^{\hat{\delta}}_z$ of a Bessel process, is given by:
\begin{equation} \label{loiT0}
   \widehat{\varphi}_{\frac{\mu^2}{8}\downarrow}(z) := E^{\hat{\delta}}_z\left(\exp(- \frac{\mu^2}{8}T_0) \right) = 2^{1-\frac{\alpha}{2}} \Gamma(\frac{\alpha}{2})^{-1} \left( \frac{\mu z}{2}\right)^{\frac{\alpha}{2}}
K_{\frac{\alpha}{2}} \left( \frac{\mu z}{2} \right). 
\end{equation}
Using \eqref{eq2.12}, we can recover the expression of $\varphi_{\theta\downarrow}$ for $\dis\theta=\frac{\delta\mu}{2}$ obtained in \eqref{eqlaplace} using the identity \eqref{eq2.5}.
\sub We now develop a discussion
similar to that made in \eqref{eq2.7}, but with the downwards arrows
$\downarrow$ now changed into upwards arrows $\uparrow$ (for the
definition of these pushed upwards and downwards processes
obtained from a
diffusion, see Pitman-Yor \cite{PY2}).\\
The analogue of formula \eqref{eq2.9} is now :
\begin{eqnarray}\label{eq2.14}
    &&Q_{z|_{\mathcal{F}_t}}^{\delta,\,\mu;\,\theta\uparrow}\\
    &=&\frac{\varphi_{\theta\uparrow}(Z_t)}{\varphi_{\theta\uparrow}(z)}\exp(-\theta t)\exp\left(-\frac{\mu}{2}(Z_t-\delta t-z)-\frac{\mu^2}{2}\int_0^tdsZ_s\right)\textrm{\huge.}\,Q_{z|_{\mathcal{F}_t}}^\delta\nonumber
\end{eqnarray}
and we note again that, precisely for :
$\dis\theta=\frac{\delta\mu}{2}$, this relation simplifies as :
\begin{equation}\label{eq2.15}
    Q_{z|_{\mathcal{F}_t}}^{\delta,\,\mu;\,\frac{\delta\mu}{2}\uparrow}=\frac{\exp\left(-\frac{\mu}{2}Z_t\right)\varphi_{\theta\uparrow}(Z_t)}{\exp\left(-\frac{\mu}{2}z\right)\varphi_{\theta\uparrow}(z)}\exp\left(-\frac{\mu^2}{2}\int_0^tdsZ_s\right)\textrm{\huge.}\,Q_{z|_{\mathcal{F}_t}}^\delta
\end{equation}
(We note that this formula is even simpler than \eqref{eq2.10} since here
there is no local time contribution).\\
We now continue to develop an analogous discussion to that made in
subsection \eqref{eq2.7}.\\
Thus, we time-change both sides of the absolute continuity
relation \eqref{eq2.15} with the inverse of $\dis 4\int_0^t
Z_u\,du\,,\;t\geq0$, with $\widehat{Z}$, as defined from $Z$ in
\eqref{eq2.11}. We obtain :
\begin{equation*}
    \widehat{Q}_{z|_{\widehat{\mathcal{F}}_u}}^{\delta,\,\mu;\,(\frac{\delta\mu}{2})\uparrow}=\frac{\exp\left(-\frac{\mu}{2}\widehat{Z}_u\right)\varphi_{\theta\uparrow}(\widehat{Z}_u)}{\exp\left(-\frac{\mu}{2}z\right)\varphi_{\theta\uparrow}(z)}\,\exp\left(-\frac{\mu^2u}{8}\right)\textrm{\huge.}\,\widehat{Q}_{z|_{\widehat{\mathcal{F}}_u}}^{\delta}
\end{equation*}
With obvious notation, this right-hand side may be written :
\begin{equation*}
    \frac{\widehat{\varphi}_{\frac{\mu^2}{8}\uparrow}(\widehat{Z}_u)}{\widehat{\varphi}_{\frac{\mu^2}{8}\uparrow}(z)}\,\exp\left(-\frac{\mu^2u}{8}\right)\textrm{\huge.}\,\widehat{Q}_{z|_{\widehat{\mathcal{F}}_u}}^{\delta}
\end{equation*}
with :
\begin{equation} \label{WhitBes}
    \widehat{\varphi}_{\frac{\mu^2}{8}\uparrow}(z)=e^{-\frac{\mu}{2}z}\varphi_{\theta\uparrow}(z)
\end{equation}
 $\widehat{Q}_{z}^{\delta,\,\mu;\,(\frac{\delta\mu}{2})\uparrow}$ is the distribution of a Bessel process of dimension $\hat{\delta} = 2-\alpha$ with drift $\frac{\mu^2}{8}$, i.e.
 $$ \widehat{Q}_{z}^{\delta,\,\mu;\,(\frac{\delta\mu}{2})\uparrow} = P_z^{\hat{\delta}; \frac{\mu^2}{8} \uparrow}.$$
 The analytical counterpart of \eqref{WhitBes} is the companion formula of \eqref{eq2.5} (see Appendix):
 $$  M_{0,\,-\frac{\alpha}{2}}(\xi)= 4^{-\frac{\alpha}{2}} \Gamma(1-\frac{\alpha}{2} ) \sqrt{\xi}\;I_{-\frac{\alpha}{2}}\left(\frac{\xi}{2}\right); $$
 while the companion formula of \eqref{loiT0} is:
 $$  \widehat{\varphi}_{\frac{\mu^2}{8}\uparrow}(z) := \frac{1}{E_0^{\hat{\delta}}\left(\exp(- \frac{\mu^2}{8}T_z) \right)}  = 2^{-\frac{\alpha}{2}} \Gamma(1- \frac{\alpha}{2}) \left( \frac{\mu z}{2}\right)^{\frac{\alpha}{2}}
I_{-\frac{\alpha}{2}} \left( \frac{\mu z}{2} \right) $$
a well-known formula which goes back to Kent \cite{Ke}, Pitman-Yor \cite{PY2}.
\section*{Appendix : On the Whittaker and Bessel - \\Mc Donald
functions} The following formulae involving these classical
special functions are found in Lebedev  \cite{L}, to which we refer
with numberings such as : $(N)_*$ ...
\begin{enumerate}
    \item [a)] The Whittaker functions $M_{k,\,\mu}(z)$ and $W_{k,\,\mu}(z)$ are a pair of
    solutions of Whittaker's equation :
    \begin{equation*}
        u''+\left(-\frac{1}{4}+\frac{k}{z}+\frac{(\frac{1}{4}-\mu^2)}{z^2}\right)u=0
    \end{equation*}
    (p. 279$_*$).
    \item [b)] $W_{k,\,\mu}$ admits the integral representation :
    \begin{equation*}
        W_{k,\,\mu}(z)=\frac{z^k\,e^{-\frac{z}{2}}}{\Gamma(\mu-k+\frac{1}{2})}\int_0^\infty
        e^{-t}\,t^{\mu-k-\frac{1}{2}}\left(1+\frac{t}{z}\right)^{\mu-k+\frac{1}{2}}\ dt
    \end{equation*}
    (see Problem 17$_*$, p. 279$_*$).
    \item [c)] \begin{equation*}
        W_{0,\,\mu}(z)=\sqrt{\frac{z}{\pi}}\;K_\mu\left(\frac{z}{2}\right)
    \end{equation*}
    (see Problem 19$_*$, p. 279$_*$).
    \item [d)] In terms of the confluent hypergeometric function
    $\Psi$, there are the relations :
    \begin{equation*}
        W_{k,\,\mu}(z)=z^{\mu+\frac{1}{2}}\,e^{-\frac{z}{2}}\Psi\left(\frac{1}{2}-k+\mu,\;2\mu+1;\;z\right)
    \end{equation*}
    (see (9.13.16)$_*$, p. 274$_*$).
    \begin{equation*}
        K_\mu(z)=\sqrt{\pi}(2z)^\mu\,e^{-z}\Psi\left(\mu+\frac{1}{2},\;2\mu+1;\;2z\right)
    \end{equation*}
    (see (9.13.15)$_*$, p. 274$_*$).\\
    Taking $k=0$ in the above formula for $W_{k,\,\mu}$, one recovers c).
 \item[e)] In terms of the confluent hypergeometric function
    $\Phi$, there are the relations :
    \begin{equation*}
        M_{k,\,\mu}(z)=z^{\mu+\frac{1}{2}}\,e^{-\frac{z}{2}}\Phi\left(\frac{1}{2}-k+\mu,\;2\mu+1;\;z\right)
    \end{equation*}
    (see (9.13.16)$_*$, p. 274$_*$).
    \begin{equation*}
        I_\mu(z)=\frac{(z/2)^\mu}{ \Gamma(\mu +1)} \,e^{-z}\Phi\left(\mu+\frac{1}{2},\;2\mu+1;\;2z\right)
    \end{equation*}
    (see (9.13.14)$_*$, p. 274$_*$).\\
    Taking $k=0$ in the above formula for $M_{k, \mu}$, we obtain:
 $$  M_{0,\,\mu}(z)= 4^{\mu} \Gamma(\mu +1) \sqrt{z}\;I_{\mu}\left(\frac{z}{2}\right) $$ 
\end{enumerate}

\end{document}